\documentclass[11pt]{article}
\usepackage{epsf}

\usepackage{amsthm}
\usepackage{amssymb}
\usepackage{amsmath}

\usepackage[top=1.5cm,bottom=1.5cm,left=2.5cm,right=2.5cm,includehead,includefoot
           ]{geometry}

\makeatletter

\numberwithin{equation}{section}
\newtheorem{theorem}{Theorem}[section]
\newtheorem{lemma}[theorem]{Lemma}
\newtheorem{corollary}[theorem]{Corollary}

\title{ Mappings on some reflexive algebras characterized by action on zero products or Jordan zero products }
\author{\begin{tabular}{c}Yunhe Chen and Jiankui Li
\footnote{Corresponding author.
E-mail address:jiankuili@yahoo.com }\\
{\small\it Department of Mathematics, East China University of
Science and Technology}\\
{\small\it Shanghai 200237, P. R. China}
\end{tabular}}
\date{}
\begin{document}
\maketitle \abstract Let
 $\mathcal{L}$ be  a subspace lattice
 on a Banach space $X$ and let $\delta:\mathrm{Alg}\mathcal{L}\rightarrow B(X)$ be a linear
 mapping.
If $\vee\{L\in \mathcal{L}: L_-\nsupseteq L\}=X$ or
$\wedge\{L_-:L\in \mathcal{L}, L_-\nsupseteq L\}=(0)$, we show that
the following three conditions are equivalent: (1)
$\delta(AB)=\delta(A)B+A\delta(B)$ whenever  $AB=0$; (2)
$\delta(AB+BA)=\delta(A)B+A\delta(B)+\delta(B)A+B\delta(A)$ whenever
$AB+BA=0$; (3) $\delta$ is a generalized derivation and
$\delta(I)\in (\mathrm{Alg}\mathcal{L})^\prime$. If
  $\vee\{L\in \mathcal{L}: L_-\nsupseteq L\}=X$ or
$\wedge\{L_-:L\in \mathcal{L}, L_-\nsupseteq L\}=(0)$ and  $\delta$
satisfies
$\delta(AB+BA)=\delta(A)B+A\delta(B)+\delta(B)A+B\delta(A)$ whenever
$AB=0$,
 we obtain that  $\delta$ is a generalized derivation and
$\delta(I)A\in(\mathrm{Alg}\mathcal{L})^\prime$ for every $A\in
\mathrm{Alg}\mathcal{L}$. We also prove that if $\vee\{L\in
\mathcal{L}: L_-\nsupseteq L\}=X$ and $\wedge\{L_-:L\in \mathcal{L},
L_-\nsupseteq L\}=(0)$, then $\delta$ is a local generalized
derivation if and only if $\delta$ is a generalized derivation.

\


\


\

\section{Introduction}
Throughout this paper, let $X$ be a Banach space over the real or
complex  field $\mathbb{F}$ and $X^*$ be the topological dual of
$X$.  When $X$ is a Hilbert space, we change it to $H$. We denote by
$B(X)$ the set of all bounded linear operators on $X$. For $A\in
B(X)$, we denote by $A^*$ the adjoint of $A$. A \emph{subspace} of
$X$ means a norm closed linear manifold. For a subset $L\subseteq
X$, denote by $L^\perp$ the annihilator of $L,$ that is,
$L^\perp=\{f\in X^*:f(x)=0 ~~\mathrm{for}~~ \mathrm{all}~~ x\in
L\}$. By a \emph{subspace lattice} on $X$, we mean a collection
$\mathcal{L}$ of
 subspaces of $X$ with $(0)$ and $X$ in $\mathcal{L}$
such that for every family $\{M_r\}$ of elements of $\mathcal{L}$,
both $\wedge M_r$ and $\vee M_r$ belong to $\mathcal{L}$, where
$\wedge M_r$ denotes the intersection of $\{M_r\}$ and $\vee M_r$
denotes the closed linear span of $\{M_r\}$. We use
$\mathrm{Alg}\mathcal{L}$ to denote the algebra of operators in
$B(X)$ that leave members of $\mathcal{L}$ invariant.

Let $x\in X$ and $f\in X^*$  be non-zero. The rank-one operator
$x\otimes f$ is defined by $y\mapsto f(y)x$ for $y\in X.$ If
$\mathcal{L}$ is a subspace lattice on $X$ and $E\in\mathcal{L}$, we
define
$$E_-=\vee\{F\in \mathcal{L}:F\nsupseteq E\}~,~E_+=\wedge\{F\in \mathcal{L}:F\nsubseteq E\}$$
and
$$\mathcal{J}_\mathcal{L}=\{L\in \mathcal{L}:L\neq (0)~~\mathrm{and}~~L_-\neq X\}
~,~\mathcal{P}_\mathcal{L}=\{L\in \mathcal{L}:L_-\nsupseteq L\}.$$
It is obvious that $\mathcal{P}_\mathcal{L}\subseteq
\mathcal{J}_\mathcal{L}.$ It is well known that a rank one operator
$x\otimes f\in \mathrm{Alg}\mathcal{L}$ if and only if there exists
a $K\in \mathcal{J}_\mathcal{L}$ such that $x\in K$ and $f\in
K_-^\perp$.  A subspace lattice $\mathcal{L}$ is called a
\emph{completely distributive lattice} if $L=\vee\{E\in
\mathcal{L}:E_-\nsupseteq L\}$ for every $L\in\mathcal{L}$ (see
\cite{Strongly ref});  $\mathcal{L}$ is called a
\emph{$\mathcal{J}$-subspace lattice} if $L\wedge L_-=(0)$ for every
$L\in\mathcal{J}_\mathcal{L}$, $X=\vee\{L:L\in
\mathcal{J}_\mathcal{L}\}$ and $\wedge\{L_-:L\in
\mathcal{J}_\mathcal{L}\}=(0)$ (see \cite{W. Longstaff O. Panaia}).
A totally ordered subspace lattice $\mathcal{N}$ is called a
\emph{nest}.   Recall that  $\mathcal{N}$ is a \emph{discrete nest}
if a nest $\mathcal{N}$ satisfies $N_-\neq N$ for every non-trivial
subspace $N$ in $\mathcal{N}$.

We say that $\mathcal{L}$ is a
\emph{$\mathcal{P}$-subspace lattice} on $X$  if $\vee\{L:L\in
\mathcal{P}_\mathcal{L}\}=X$ or $\wedge\{L_-:L\in
\mathcal{P}_\mathcal{L}\}=(0)$. It is obvious that this class of subspace lattices contains
$\mathcal{J}$-subspace lattices, discrete nests and subspace
lattices with $X_-\neq X$ or $(0)_+\neq (0).$  The  following
example is also a $\mathcal{P}$-subspace lattice.

{\flushleft{\textbf{Example 1.1.}} Let $\{e_n:n\in\mathbb{N}\}$ be
an orthonormal basis of $H$, $P_n= \mathrm{span}\{e_i: i=1,...,n\}$,
$\xi=\sum_{n=1}^\infty\frac{1}{n}e_n$ and $P_\xi$ be the orthogonal
projection from $H$ onto the one-dimensional subspace of $H$
generated by $\xi$.}{ It follows from \cite[Theorem 2.11]{L. Wang W.
Yuan} and \cite[Lemma 3.2]{C. Hou} that $\mathcal{L}=\{0,I,P_n,P_\xi,P_\xi\vee
P_n:n=1,2,\cdots\}$ is a reflexive
$\mathcal{P}$-subspace lattice.}

In a Hilbert space, we disregard the distinction between a closed
subspace and the orthogonal projection onto it. A subspace lattice
on a Hilbert space $H$ is called a \emph{commutative subspace
lattice} (or \emph{CSL} for short) if it consists of mutually
commuting projections. In the paper, we assume that $H$ is a complex
separable Hilbert space.

Let $\delta$ be a linear mapping from a unital algebra $\mathcal{A}$
into an $\mathcal{A}$-bimodule $\mathcal{M}$. Recall that $\delta$
is a \emph{derivation} (respectively  \emph{generalized derivation}) if
$\delta(AB)=\delta(A)B+A\delta(B)$ (respectively
$\delta(AB)=\delta(A)B+A\delta(B)-A\delta(I)B$) for all $A,B$ in
$\mathcal{A}$. We say that $\delta$ is \emph{derivable} at $Z\in
\mathcal{A}$ if $\delta(AB)=\delta(A)B+A\delta(B)$ for any $A,B\in
\mathcal{A}$ with $AB=Z$; $\delta$ is \emph{Jordan derivable} at
$Z\in \mathcal{A}$ if
$\delta(AB+BA)=\delta(A)B+A\delta(B)+\delta(B)A+B\delta(A)$ for any
$A,B\in \mathcal{A}$ with $AB+BA=Z$.
If $\delta(AB+BA)=\delta(A)B+A\delta(B)+\delta(B)A+B\delta(A)$ for any
$A,B\in \mathcal{A}$ with $AB=0$, we say that $\delta$ has
\emph{WJD} (weak Jordan derivation) property.

In recent years, there have been a number of papers on the study of
conditions under which derivations and Jordan derivations of
operator algebras can be completely determined by the action on some
subsets of operator algebras (for example, see \cite{R. An J. Hou,
M. Chebotar, Meiyan Jiao Jinchuan Hou, W. Jing S. Lu P. Li, S. Zhao
J. Zhu}). For instance, Zhao and Zhu in \cite{S. Zhao J. Zhu} showed
that every linear mapping $\delta$ from a triangular algebra
$\mathcal{T}$ into itself satisfying WJD property is a derivation.
In \cite{Meiyan Jiao Jinchuan Hou}, Jiao and Hou proved that every
additive mapping $\delta$ derivable or Jordan derivable at zero
point on some nest algebras has the form $\delta(A)=\tau(A)+cA$ for
some additive derivation $\tau$ and some scalar $c\in\mathbb{F}$.

The purpose of this paper is to consider some mappings which behave like
derivations on $\mathcal{P}$-subspace lattice algebras and
completely distributive commutative subspace lattice (CDCSL)
algebras.

In Section 2, we show that every linear (respectively bounded linear) mapping
$\delta$ on $\mathcal{P}$-subspace lattice (respectively CDCSL) algebras
Jordan derivable at zero point is a generalized derivation and
$\delta(I)\in(\mathrm{Alg}\mathcal{L})^\prime$.

 In Section 3, for a $\mathcal{P}$-subspace lattice algebra
$\mathrm{Alg}\mathcal{L}$, we obtain that  $\delta$ satisfies WJD
property if and only if $\delta$ is a generalized derivation and
$\delta(I)A\in(\mathrm{Alg}\mathcal{L})^\prime$ for every $A\in
\mathrm{Alg}\mathcal{L}$.

    In Section 4, we  investigate derivable mappings at zero point and
 some linear mappings which behave like
 left (respectively right) multipliers, isomorphisms or local generalized
derivations on $\mathcal{P}$-subspace lattice algebras. One of the
main results of the section is that if $\vee\{L\in \mathcal{L}:
L_-\nsupseteq L\}=X$ and $\wedge\{L_-:L\in \mathcal{L},
L_-\nsupseteq L\}=(0)$, then $\delta$ is a local generalized
derivation from $\mathrm{Alg}\mathcal{L}$ into $B(X)$ if and only if
$\delta$ is a generalized derivation.

The following proposition will be used in our proofs.

{\flushleft{\textbf{Proposition 1.2}} (\cite [Proposition 1.1]{F. Lu B. Liu})\textbf{.}
\emph{Let $E$ and $F$ be non-zero subspaces of $X$ and $X^*$, respectively.}}
{\emph{Let $\Phi:E\times F\rightarrow B(X)$ be a bilinear mapping such that
 $\Phi(x,f)ker(f)\subseteq \mathbb{F}x$ for all $x\in E$ and $f\in F$.}
\emph{Then there exist two linear mappings $T:E\rightarrow X$ and
$S:F\rightarrow X^*$ such that $\Phi(x,f)=Tx\otimes f+x\otimes Sf$
for all $x\in E$ and $f\in F$.}
}

\section{Jordan derivable Mappings at zero point}

\

The following lemma is included in the proof of \cite[Theorem 3.1]{Meiyan Jiao Jinchuan Hou}.
We leave the proof to readers.
\begin{lemma}\label{2001}
If $\delta$ is Jordan derivable at zero point from a unital algebra $\mathcal{A}$ into its unital bimodule,
then for any idempotents $P$ and $Q$ in $\mathcal{A}$, the following hold:\\
$(1)~~\delta(I)P=P\delta(I);$\\
$(2)~~\delta(P)=\delta(P)P+P\delta(P)-P\delta(I);$\\
$(3)~~\delta(PQ+QP)=\delta(P)Q+P\delta(Q)+\delta(Q)P+Q\delta(P)-\delta(I)(PQ+QP).$
\end{lemma}

For a subspace lattice $\mathcal{L}$ and a subspace $E\in \mathcal{P}_\mathcal{L}$,
we denote by $\mathcal{T}_E$ the ideal span$\{x\otimes f:x\in E,f\in E_-^\perp\}$ of $\mathrm{Alg}\mathcal{L}$.

\begin{lemma}\label{2002}
If $\mathcal{L}$ is a subspace lattice on $X$ and $E$ is in $\mathcal{P}_\mathcal{L}$,
 then for every $x$ in $E$ and  every $f$ in $E_-^{\perp}$,
$x\otimes f$  is a linear combination of idempotents in $\mathcal{T}_E$.
\end{lemma}
\begin{proof}[\bf Proof]
Suppose $f(x)\neq 0$, then $x\otimes f=f(x)(\frac{1}{f(x)}x\otimes f),$
where $\frac{1}{f(x)}x\otimes f$ is an idempotent in $\mathcal{T}_E$.
Suppose $f(x)=0$.
Since $E\in\mathcal{P}_\mathcal{L},$ there exist $z\in E$ and $g\in E_-^\perp$ such that $g(z)=1.$ \

Case 1. If $g(x)=\mu\neq 0$, then $x\otimes f=x\otimes(\frac{1}{\mu}g+f)-x\otimes\frac{1}{\mu}g,$
where $x\otimes(\frac{1}{\mu}g+f)$ and $x\otimes\frac{1}{\mu}g$ are idempotents in $\mathcal{T}_E$.

Case 2. If $f(z)=\lambda\neq 0$, then $x\otimes f=(x+\frac{1}{\lambda}z)\otimes f-\frac{1}{\lambda}z\otimes f,$
where $(x+\frac{1}{\lambda}z)\otimes f$ and $\frac{1}{\lambda}z\otimes f$ are idempotents in $\mathcal{T}_E$.

Case 3. If $f(z)=g(x)=0$, then $x\otimes f=\frac{1}{4}((z+x)
\otimes(g+f)+(z-x)\otimes(g-f)-(z+x)\otimes(g-f)-(z-x)\otimes(g+f)),$
where $(z+x)\otimes(g+f)$, $(z-x)\otimes(g-f)$, $(z+x)\otimes(g-f)$
 and $(z-x)\otimes(g+f)$ are idempotents in $\mathcal{T}_E$. The proof is complete.
\end{proof}

\begin{lemma}\label{2003}
Let $\mathcal{L}$ be a subspace lattice on $X$, $E$ be in $P_\mathcal{L}$ and $\delta$ be
a linear mapping from $\mathrm{Alg}\mathcal{L}$ into $B(X)$. If $\delta$ is Jordan derivable at zero point,
then for every idempotent P in $\mathcal{T}_E$ and every $A$ in $\mathrm{Alg}\mathcal{L}$, the following hold:\\
$(1)~~\delta(AP+PA)=\delta(A)P+A\delta(P)+\delta(P)A+P\delta(A)-\delta(I)(AP+PA)$;\\
$(2)~~\delta(PAP)=\delta(P)AP+P\delta(A)P+PA\delta(P)-2\delta(I)PAP.$
\end{lemma}

\begin{proof}[\bf Proof]
(1)~~For every idempotent $P\in\mathcal{T}_E$ and every $A\in \mathrm{Alg}\mathcal{L}$, since
$P^\perp AP^\perp P+PP^\perp AP^\perp=0$, by assumption we have
\begin{eqnarray*}
\delta(P^\perp AP^\perp)P+P^\perp AP^\perp\delta(P)+\delta(P)P^\perp AP^\perp+P\delta(P^\perp AP^\perp)=0.
\end{eqnarray*}
Since $A-P^\perp AP^\perp=PA+P^\perp AP\in \mathcal{T}_E$, it follows from Lemmas \ref{2001} and \ref{2002} that
\begin{eqnarray*}
\delta(AP+PA)&=&\delta((A-P^\perp AP^\perp)P+P(A-P^\perp AP^\perp))\\
&=&\delta(A-P^\perp AP^\perp)P+(A-P^\perp AP^\perp)\delta(P)+\delta(P)(A-P^\perp AP^\perp)\\
&&+P\delta(A-P^\perp AP^\perp)-\delta(I)(AP+PA)\\
&=&\delta(A)P+A\delta(P)+\delta(P)A+P\delta(A)-\delta(I)(AP+PA)\\
&&-(\delta(P^\perp AP^\perp)P+P^\perp AP^\perp\delta(P)+\delta(P)P^\perp AP^\perp+P\delta(P^\perp AP^\perp))\\
&=&\delta(A)P+A\delta(P)+\delta(P)A+P\delta(A)-\delta(I)(AP+PA).
\end{eqnarray*}
(2) The substitution $AP+PA$ for $A$ in (1)  gives (2).
\end{proof}

One of the main results of this section is the following theorem.

\begin{theorem}\label{2004}
Let $\mathcal{L}$ be a subspace lattice on $X$ such that
$\vee\{L:L\in \mathcal{P}_\mathcal{L}\}=X$ and $\delta$  be a linear
mapping from $\mathrm{Alg}\mathcal{L}$ into $B(X)$. Then  $\delta$
is Jordan derivable at zero point if and only if $\delta$ is a
generalized derivation and
$\delta(I)\in(\mathrm{Alg}\mathcal{L})^\prime$, where
$(\mathrm{Alg}\mathcal{L})^\prime$ is the commutant of
$\mathrm{Alg}\mathcal{L}$ in $B(X)$. In particular, if
$\delta(I)=0,$ then $\delta$ is Jordan derivable at zero point if
and only if $\delta$ is a derivation.
\end{theorem}
\begin{proof}[\bf Proof]
The sufficiency is obvious, so we only need to prove the necessity.
Let $E\in\mathcal{P}_\mathcal{L}$, $z\in E$ and $g\in E_-^\perp$
with $g(z)=1$.
We divide the proof into several claims.\\
\textbf{Claim 1.} $\delta(I)\in(\mathrm{Alg}\mathcal{L})^\prime.$

For all $x\in E$, $f\in E_-^\perp$ and
$T\in\mathrm{Alg}\mathcal{L}$, by Lemmas \ref{2001} and \ref{2002},
we have $\delta(I)Tx\otimes f=Tx\otimes
f\delta(I)=T\delta(I)x\otimes f.$ That is, $\delta(I)Tx=T\delta(I)x$
for every  $x\in E$. Since $\vee\{E:  E\in
\mathcal{P}_\mathcal{L}\}=X$, it follows that $\delta(I)\in
(\mathrm{Alg}\mathcal{L})^\prime.$

Now define $\tau(A)=\delta(A)-\delta(I)A$ for $A\in \mathrm{Alg}\mathcal{L}$. It is easy to see that
$\tau$ is Jordan derivable at zero point and $\tau(I)=0$.\\
\textbf{Claim 2.} $\tau(x\otimes f)ker(f)\subseteq \mathbb{F}x,$ for all  $x\in E$ and $f\in E_-^\perp.$

Case 1. If $f(x)=\mu\neq0$, then by Lemma \ref{2001}, we have
\begin{eqnarray*}
\tau(\frac{1}{\mu}x\otimes f))=\tau(\frac{1}{\mu}x\otimes f)
(\frac{1}{\mu}x\otimes f)+(\frac{1}{\mu}x\otimes f)\tau(\frac{1}{\mu}x\otimes f).
\end{eqnarray*}
Thus $\tau(x\otimes f)ker(f)\subseteq \mathbb{F}x.$

Case 2. If $f(x)=0$ and $f(z)\neq0$, then by Case 1, for every $y\in ker(f)$, we have
\begin{eqnarray*}
\tau((z+x)\otimes f)y&=&\lambda_1(z+x),\\
\tau((z-x)\otimes f)y&=&\lambda_2(z-x),\\
\tau(z\otimes f)y&=&\lambda_3 z,
\end{eqnarray*}
for some $\lambda_1$, $\lambda_2$ and $\lambda_3\in \mathbb{F}$. By
the above equations, it follows that
$$2\lambda_3z=(\lambda_1+\lambda_2)z+(\lambda_1-\lambda_2)x,$$
and the independence of $z$ and $x$ implies $\lambda_1=\lambda_2=\lambda_3$.
Hence
$$\tau(x\otimes f)y=\tau((z+x)\otimes f)y-\tau(z\otimes f)y=\lambda_1 x.$$
This means $\tau(x\otimes f)ker(f)\subseteq \mathbb{F}x.$\

Case 3. Suppose that $f(x)=0$ and $f(z)=0.$  Since $z\otimes (g+f)$
and $z\otimes(g-f)$ are idempotents in $\mathcal{T}_E$, it follows
from Lemma \ref{2003} that
\begin{eqnarray*}
&&\tau((z\otimes (g+f))(x\otimes g)(z\otimes(g+f)))\\
&&~~~~~~~~~~~~~=\tau(z\otimes(g+f))(x\otimes g)(z\otimes(g+f))+(z\otimes (g+f))\tau(x\otimes g)(z\otimes(g+f))\\
&&~~~~~~~~~~~~~~~~+(z\otimes (g+f))(x\otimes g)\tau(z\otimes(g+f)),
\end{eqnarray*}
\begin{eqnarray*}
&&\tau((z\otimes (g-f))(x\otimes g)(z\otimes(g-f)))\\
&&~~~~~~~~~~~~~=\tau(z\otimes(g-f))(x\otimes g)(z\otimes(g-f))+(z\otimes (g-f))\tau(x\otimes g)(z\otimes(g-f))\\
&&~~~~~~~~~~~~~~~~+(z\otimes (g-f))(x\otimes g)\tau(z\otimes(g-f)),
\end{eqnarray*}
and
\begin{eqnarray*}
&&\tau((z\otimes g)(x\otimes g)(z\otimes g))\\
&&~~~~~~~~~~~=\tau(z\otimes g)(x\otimes g)(z\otimes g)+(z\otimes g)\tau(x\otimes g)
(z\otimes g)+(z\otimes g)(x\otimes g)\tau(z\otimes g).
\end{eqnarray*}
From the above three equations, we have
\begin{eqnarray*}
0&=&\tau((z\otimes f)(x\otimes g)(z\otimes f))\\
&=&\tau(z\otimes f)(x\otimes g)(z\otimes f)+(z\otimes f)\tau(x\otimes g)(z\otimes f)
+(z\otimes f)(x\otimes g)\tau(z\otimes f)\\
&=&\tau(z\otimes f)(x\otimes f)+(z\otimes f)\tau(x\otimes g)(z\otimes f).
\end{eqnarray*}
Thus
\begin{eqnarray}
\tau(z\otimes f)x=-f(\tau(x\otimes g)z)z.\label{201}
\end{eqnarray}
Hence by (\ref{201}), Lemmas 2.2 and 2.3, it follows that
\begin{eqnarray*}
&&\tau(x\otimes f)=\tau((z\otimes f)(x\otimes g)+(x\otimes g)(z \otimes f))\\
&&~~~~~~~~~~~=-f(\tau(x\otimes g)z)z\otimes g+(z\otimes f)\tau(x\otimes g)\\
&&~~~~~~~~~~~~~~+\tau(x\otimes g)(z\otimes f)+(x\otimes g)\tau(z\otimes f).
\end{eqnarray*}
Let $y$ be in $ker(f)$. Applying the above equations to $y$ gives
\begin{eqnarray}
\tau(x\otimes f)y=-g(y)f(\tau(x\otimes g)z)z+f(\tau(x\otimes g)y)z+g(\tau(z\otimes f)y)x.\label{202}
\end{eqnarray}
Notice that (\ref{202}) is valid for all $z\in E$ satisfying $g(z)=1$ and $f(z)=0$.
If $g(x)=\mu\neq 0,$ replacing $z$ by $\frac{1}{\mu}x$ in (\ref{202}),
we have $\tau(x\otimes f)y\in \mathbb{F}x$.
If $g(x)=0,$ by the proof of \cite[Lemma 2.3]{Lu F. reflexive}, we have
$g(y)f(\tau(x\otimes g)z)-f(\tau(x\otimes g)y)=0,$
whence $\tau(x\otimes f)y=g(\tau(z\otimes f)y)x\in \mathbb{F}x$.\\
\textbf{Claim 3.} $\tau$ is a derivation.

By Claim 2 and Proposition 1.2, there exist linear mappings $T:E\rightarrow X$
and $S:E_-^\perp\rightarrow X^*$ such that
\begin{eqnarray}
\tau(x\otimes f)=Tx\otimes f+x\otimes Sf,\label{203}
\end{eqnarray}
for all $x\in E$ and $f\in E_-^\perp.$ It follows from Lemmas
\ref{2002} and \ref{2003} that for every $A\in
\mathrm{Alg}\mathcal{L}$,
\begin{eqnarray}
\tau(Ax\otimes g+x\otimes gA)=\tau(A)x\otimes g+A
\tau(x\otimes g)+\tau(x\otimes g)A+x\otimes g\tau(A).\label{204}
\end{eqnarray}
By (\ref{203}) and (\ref{204}), we have
\begin{eqnarray*}
&&TAx\otimes g+Ax\otimes Sg+Tx\otimes A^*g+x\otimes SA^*g\\
&&~~~~~~~~~~~~=\tau(A)x\otimes g+ATx\otimes g+Ax\otimes Sg+Tx\otimes A^*g+x\otimes A^*Sg+x\otimes\tau(A)^*g.
\end{eqnarray*}
That is,
$$(\tau(A)+AT-TA)x\otimes g=x\otimes(SA^*-\tau(A)^*-A^*S)g.$$
Thus there exists a linear mapping
$\lambda:\mathrm{Alg}\mathcal{L}\rightarrow \mathbb{F}$ such that
\begin{eqnarray}
\tau(A)x=(TA-AT)x+\lambda(A)x,\label{205}
\end{eqnarray}
for all $A\in \mathrm{Alg}\mathcal{L}$ and $x\in E$.
Hence by (\ref{205}), for all $A,B$ in $\mathrm{Alg}\mathcal{L}$ and $x$ in $E$,
\begin{eqnarray}
\tau(AB)x=(\tau(A)B+A\tau(B))x+\lambda(AB)x-\lambda(A)Bx-\lambda(B)Ax.\label{206}
\end{eqnarray}

In the following, we show $\lambda(A)=0$ for every $A\in \mathrm{Alg}\mathcal{L}.$
Putting $A=B=z\otimes g$ and $x=z$ in (\ref{206}) gives
$\lambda(z\otimes g)=g(\tau(z\otimes g)z),$
and  Lemma \ref{2001} (2) implies $g(\tau(z\otimes g)z)=0.$
Hence
\begin{eqnarray}
\lambda(z\otimes g)=0.\label{207}
\end{eqnarray}
Notice that (\ref{207}) is valid for all $z$ in $E$ and $g$ in $E_-^\perp$ satisfying
 $g(z)=1$. Now fix $z\in E$ and $g\in E_-^\perp$ such that $g(z)=1$.
Thus for all $f\in E_-^\perp$, if $f(z)=\mu\neq0$, then $\lambda(z\otimes f)$=
$\mu\lambda( z\otimes \frac{1}{\mu}f)=0;$ if  $f(z)=0$,
then $\lambda(z\otimes f)=\lambda(z\otimes(g+f))-\lambda(z\otimes g)=0.$ Hence
$\lambda(z\otimes f)=0$ for every $f\in E_-^\perp.$ Similarly,
we have $\lambda(x\otimes g)=0$ for every $x\in E.$ Now for every $A\in \mathrm{Alg}
\mathcal{L}$, by (\ref{206}), we have
\begin{eqnarray}
\tau(Az\otimes g)z=\tau(A)z+A\tau(z\otimes g)z-\lambda(A)z \label{208}
\end{eqnarray}
and
\begin{eqnarray}
\tau(z\otimes gA)z=\tau(z\otimes g)Az+g(\tau(A)z)z-\lambda(A)z.\label{209}
\end{eqnarray}
By Lemma \ref{2003} (1), we have
\begin{eqnarray}
\tau(Az\otimes g+z\otimes gA)z=\tau(A)z+A\tau(z\otimes g)z+\tau(z\otimes g)Az+g(\tau(A)z)z.\label{210}
\end{eqnarray}
Combining (\ref{208}), (\ref{209}) and (\ref{210}) gives $\lambda(A)=0$ for every $A\in \mathrm{Alg}\mathcal{L}$.
 Then by (\ref{206}), we obtain
$$\tau(AB)x=(\tau(A)B+A\tau(B))x,$$
for all $A,B\in\mathrm{ Alg}\mathcal{L}$ and $x\in E.$
Since $\vee\{L:L\in \mathcal{P}_\mathcal{L}\}=X$, it follows that $\tau$ is a derivation.
By $\delta(A)=\tau(A)+\delta(I)A$, it is easy to show  that $\delta$ is a generalized derivation.
\end{proof}

Applying the ideas in the proof of Theorem \ref{2004}, we can obtain
the following result.

\begin{theorem}\label{2005}
Let $\mathcal{L}$ be a subspace lattice on $X$ such that $\wedge\{L_-:L\in \mathcal{P}_\mathcal{L}\}=(0)$
and $\delta$ be a linear mapping from $\mathrm{Alg}\mathcal{L}$ into $B(X)$.
Then  $\delta$ is Jordan derivable at zero point if and only if $\delta$ is a generalized
derivation and $\delta(I)\in(\mathrm{Alg}\mathcal{L})^\prime.$
In particular, if $\delta(I)=0,$ then $\delta$ is Jordan derivable at zero point if and only
if $\delta$ is a derivation.
\end{theorem}
\begin{proof}[\bf Proof]
We only prove the necessity.
Let $x\mapsto \hat{x}$ be the canonical mapping from $X$ into $X^{**}$, then $(x\otimes f)^*=f\otimes \hat{x}$
for all $x\in X$ and $f\in X^*$.
The hypothesis $\wedge\{L_-:L\in \mathcal{P}_\mathcal{L}\}=(0)$ implies that $\vee\{L_-^\perp:L\in
\mathcal{P}_\mathcal{L}\}=X^*$.
With a proof similar to the proof of Theorem \ref{2004},
we have $\delta(I)\in (\mathrm{Alg}\mathcal{L})^\prime$.
Let $\tau(A)=\delta(A)-\delta(I)A$ for $A\in \mathrm{Alg}\mathcal{L}$. Then  $\tau$ is Jordan
derivable at zero point and $\tau(I)=0.$
In the following, we show $\tau$ is a derivation.
Let $E\in \mathcal{P}_\mathcal{L}$. We choose $z\in E$ and $g\in E_-^\perp$ such that $g(z)=1.$ One
can easily verify that for all $x\in E$
and $f\in E_-^\perp$,
$\tau(x\otimes f)^*ker(\hat{x})\subseteq \mathbb{F}f.$ Let $\Phi(f,\hat{x})=\tau(x\otimes f)^*$
for all $x\in E$ and $f\in E_-^\perp$.
Then $\Phi$ is a bilinear mapping from $E_-^\perp\times \hat{E}$ into $B(X^*)$, where $\hat{E}=\{\hat{x}:x\in E\}$.
Hence there exist linear mappings $T:E_-^\perp\rightarrow X^*$ and $S: \hat{E}\rightarrow X^{**}$ such that
\begin{eqnarray*}
\tau(x\otimes f)^*=\Phi(f,\hat{x})=Tf\otimes \hat{x}+f\otimes S\hat{x},
\end{eqnarray*}
for all $x\in E$ and $f\in E_-^\perp$.
Hence for $A\in \mathrm{Alg}\mathcal{L}$ and $f\in E_-^\perp$, we have that
$$(\tau(A)^*+A^*T-TA^*)f\otimes \hat{z}=f\otimes (S\widehat{Az}-\widehat{\delta(A)z}-A^{**}S\hat{z}).$$
It follows that $\tau(A)^*f=(TA^*-A^*T)f+\lambda(A)f$, where $\lambda:\mathrm{Alg}\mathcal{L}\rightarrow
 \mathbb{F}$ is a linear mapping.
Hence for all $A,B\in \mathrm{Alg}\mathcal{L}$ and $f\in E_-^\perp$,
$$\tau(AB)^*f=(B^*\tau(A)^*+\tau(B)^*A^*)f-\lambda(A)B^*f-\lambda(B)A^*f+\lambda(AB)f.$$
With a proof similar to the proof of Theorem \ref{2004},
we can prove that $\lambda(A)=0$ for every $A\in \mathrm{Alg}\mathcal{L}$.
Since $\vee\{L_-^\perp:L\in \mathcal{P}_\mathcal{L}\}=X^*$,
it follows that $\tau$ is a derivation. Hence $\delta$ is a generalized derivation.
\end{proof}

Next we investigate the bounded linear mappings which are Jordan derivable at zero point
on CDCSL algebras. Recall that a CSL algebra
$\mathrm{Alg}\mathcal{L}$ is irreducible if and only if $(\mathrm{Alg}\mathcal{L})^\prime=\mathbb{C}I,$
which is equivalent to the condition that $\mathcal{L}\cap \mathcal{L}^\bot=\{0,I\}$, where
$\mathcal{L}^\perp=\{E^\perp:E\in\mathcal{L}\}$.

\begin{lemma}[\cite{Gilfeather R. Moore}]\label{2006}
Let $\mathrm{Alg}\mathcal{L}$ be a CDCSL algebra on  $H$.
Then there exists a countable set $\{P_n:n\in\Lambda\}$ of mutually
orthogonal projections in $\mathcal{L}\cap \mathcal{L}^\bot$ such
that $\vee_n P_n=I$ and each $(\mathrm{Alg}\mathcal{L})P_n$ is an
irreducible CDCSL algebra on $P_nH$; moreover,
$\mathrm{Alg}\mathcal{L}$ can be written as a direct sum
$\mathrm{Alg}\mathcal{L}=\sum\limits_n\bigoplus(\mathrm{Alg}\mathcal{L})P_n$.
\end{lemma}

\begin{lemma}[{\cite{Lu F. CDC}}]\label{2007}
Let $\mathrm{Alg}\mathcal{L}$ be a non-trivially irreducible CDCSL
algebra  on  $H$. Then there exists a non-trivial
projection $P$ in $\mathcal{L}$ such that
$P(\mathrm{Alg}\mathcal{L})P^\bot$ is faithful, that is, for $T,S\in
\mathrm{Alg}\mathcal{L}$, $TP(\mathrm{Alg}\mathcal{L})P^\bot=\{0\}$
implies $TP=0$ and $P(\mathrm{Alg}\mathcal{L})P^\bot S=\{0\}$
implies $P^\bot S=0$.
\end{lemma}

\begin{lemma}\label{2008}
Let $\mathrm{Alg}\mathcal{L}$ be an irreducible CDCSL algebra on $H$
and let $\delta:\mathrm{Alg}\mathcal{L}\rightarrow
\mathrm{Alg}\mathcal{L}$ be a bounded linear mapping and
$\delta(I)=0$. If $\delta$ is Jordan derivable  at zero point, then
$\delta$ is a derivation.
\end{lemma}
\begin{proof}[\bf Proof]
Suppose that $\mathcal{L}$ is trivial, then $\mathrm{Alg}\mathcal{L}=B(H)$ is a von Neumann algebra.
It follows from  \cite[Theorem 3.2]{R. An J. Hou} that $\delta$ is a Jordan derivation. Since
every von Neumann algebra is a semiprime ring,
by \cite[Theorem 1]{M. Bresar}, $\delta$ is a derivation.

Suppose that $\mathcal{L}$ is non-trivial. Let $P$ be the non-trivial projection in $\mathcal{L}$
provided by Lemma \ref{2007}.
Since $P(\mathrm{Alg}\mathcal{L})P^\bot$ is faithful, by \cite[Theorem 2.1]{R. An J. Hou},
$\delta$ is a Jordan derivation.
Since every Jordan derivation on a CSL algebra is a derivation \cite[Theorem 3.2]{Lu F. CSL},
it follows that $\delta$ is a derivation.
\end{proof}

\begin{theorem}\label{2009}
Let $\mathrm{Alg}\mathcal{L}$ be a CDCSL algebra on $H$
and $\delta$ be a bounded linear mapping from $\mathrm{Alg}\mathcal{L}$ into itself. Then $\delta$ is
Jordan derivable at zero point if and only if $\delta$ is a generalized derivation and $\delta(I)
\in (\mathrm{Alg}\mathcal{L})^\prime$.
In particular, if $\delta(I)=0$, then $\delta$ is Jordan derivable at zero point if and only
 if $\delta$ is a derivation.
\end{theorem}
\begin{proof}[\bf Proof]
We only prove the necessity.
Since every rank one operator in $\mathrm{Alg}\mathcal{L}$ is a linear combination of idempotents
 in $\mathrm{Alg}\mathcal{L}$ \cite[Lemma 2.3]{D. Hadwin J. Li}
and the rank one subalgebra of $\mathrm{Alg}\mathcal{L}$ is dense in $\mathrm{Alg}\mathcal{L}$ in
 the weak topology \cite[Theorem 3]{C. Laurie W. Longstaff}, by
Lemma \ref{2001}(1), we have $\delta(I)\in (\mathrm{Alg}\mathcal{L})^\prime$.
Let $\tau(A)=\delta(A)-\delta(I)A$ for $A\in \mathrm{Alg}\mathcal{L}$. Then  $\tau$ is Jordan
derivable at zero point and $\tau(I)=0.$

Let $\mathrm{Alg}\mathcal{L}=\sum\limits_n\bigoplus(\mathrm{Alg}\mathcal{L})P_n$ be the
irreducible decomposition of $\mathrm{Alg}\mathcal{L}$ as in Lemma \ref{2006}.
Let $A$ be in $\mathrm{Alg}\mathcal{L}$ and fix an index $n$. Since $P_nAP_nP_n^\bot+P_n^\bot P_nAP_n=0,$ we have
\begin{eqnarray*}
&&~~~~~~~0=\tau(P_nAP_nP_n^\bot+P_n^\bot P_nAP_n)\\
&&~~~~~~~~~~=\tau(P_nAP_n)P_n^\bot+P_nAP_n\tau(P_n^\bot)+\tau(P_n^\bot)P_nAP_n+P_n^\bot\tau(P_nAP_n),
\end{eqnarray*}
which yields that $P_n^\bot\tau(P_nAP_n)P_n^\bot=0$. Since  $P_n\in\mathcal{L}\cap \mathcal{L}^\bot$,
there holds $\tau(AP_n)=\tau(AP_n)P_n.$ By the same
way, we obtain $\tau(AP_n^\bot)=\tau(AP_n^\bot)P_n^\bot.$ Since
\begin{eqnarray*}
0=\tau(I)=\tau(P_n+P_n^\bot)=\tau(P_n)P_n+\tau(P_n^\bot)P_n^\bot,
\end{eqnarray*}
it follows that $\tau(P_n)=0.$ Now define a linear mapping
$\tau_n:(\mathrm{Alg}\mathcal{L})P_n\rightarrow(\mathrm{Alg}\mathcal{L})P_n$ by
$$\tau_n(AP_n)=\tau(AP_n)P_n,$$ for every $A\in \mathrm{Alg}\mathcal{L}$.
It is easy to show that $\tau_n$ is bounded and Jordan derivable at zero point.
Since $(\mathrm{Alg}\mathcal{L})P_n$ is irreducible and $\tau_n(P_n)=\tau(P_n)P_n=0$,
 by Lemma \ref{2008}, $\tau_n$ is a derivation.
Hence by $\tau(A)P_n=\tau(AP_n)P_n+\tau(AP_n^\perp)P_n=\tau_n(AP_n)$, we have $\tau$
is a derivation. Thus $\delta$ is a generalized derivation.
\end{proof}

\section{ Mappings satisfying WJD property}

\

Our first result in this section says that the set of all
Jordan derivable mapping at zero point from a $\mathcal{P}$-subspace
lattice algebra into $B(X)$  is bigger than the set of all mappings
satisfying WJD property. The following lemma is included in
the proof of \cite[Lemma 2.6]{Chen}.

\begin{lemma}\label{3001}
If $\delta$ is a linear mapping satisfying WJD property from a
unital algebra $\mathcal{A}$ into its unital bimodule,
then for every idempotent $P\in \mathcal{A}$ and every $A\in \mathcal{A}$, the following hold:\\
$(1)~~\delta(I)P=P\delta(I)~~\mathrm{ and }~~\delta(P)=\delta(P)P+P\delta(P)-\delta(I)P;$\\
$(2)~~\delta(PA+AP)=\delta(P)A+P\delta(A)+\delta(A)P+A\delta(P)-\delta(I)PA-PA\delta(I);$\\
$(3)~~\delta(PA+AP)=\delta(P)A+P\delta(A)+\delta(A)P+A\delta(P)-\delta(I)AP-AP\delta(I);$\\
$(4)~~2\delta(PAP)=2\delta(P)AP+2P\delta(A)P+2PA\delta(P)-PA\delta(I)-2\delta(I)AP-AP\delta(I).$
\end{lemma}

\begin{theorem}\label{3002}
Let $\mathcal{L}$ be a subspace lattice on $X$ such that
$\vee\{L:L\in \mathcal{P}_\mathcal{L}\}=X$ and  $\delta$  be a
linear mapping from $\mathrm{Alg}\mathcal{L}$ into $B(X)$. Then
$\delta$ satisfies WJD property if and only if $\delta$ is a
generalized derivation and
$\delta(I)A\in(\mathrm{Alg}\mathcal{L})^\prime$ for every $A\in
\mathrm{Alg}\mathcal{L}$. In particular, if $\delta(I)=0$, then
$\delta$ satisfies WJD property if and only if $\delta$ is a
derivation.
\end{theorem}

\begin{proof}[\bf Proof]
Since the sufficiency is evident, we will just show the necessity.
Suppose $\delta$ satisfies WJD property. We claim that
$\delta(I)A\in (\mathrm{ Alg}\mathcal{L})^\prime$ for every $A\in
\mathrm{ Alg}\mathcal{L}$. By
 Lemma \ref{3001} (1) and
the proof of Claim 1 in Theorem \ref{2004},
we have $\delta(I)\in(\mathrm{Alg}\mathcal{L})^\prime$.
Hence by Lemma \ref{3001} (2) and (3), we have that $\delta(I)AP=PA\delta(I)$ for every idempotent
 $P\in \mathrm{Alg}\mathcal{L}$
and every $A\in\mathrm{Alg}\mathcal{L}$. Hence for all $x\in E$, $f\in E_-^\perp$ and $T\in
 \mathrm{Alg}\mathcal{L}$,
we have $\delta(I)ATx\otimes f=Tx\otimes
fA\delta(I)=T\delta(I)Ax\otimes f$. Since $\vee\{L:L\in
\mathcal{P}_\mathcal{L}\}=X$, it  follows that $\delta(I)A\in
(\mathrm{ Alg}\mathcal{L})^\prime$ for every $A\in
\mathrm{Alg}\mathcal{L}$. Let $\tau(A)=\delta(A)-\delta(I)A$ for
$A\in \mathrm{Alg}\mathcal{L}$. It is easy to show that $\tau$
satisfies WJD property and $\tau(I)=0.$ Similar to the proof of
Theorem \ref{2004}, we may show $\tau$ is a derivation and then
$\delta$
 is a generalized derivation.
\end{proof}

Similarly, we have the following theorem.

\begin{theorem}\label{3333}
Let $\mathcal{L}$ be a subspace lattice on $X$ such that
$\wedge\{L_-:L\in \mathcal{P}_\mathcal{L}\}=(0)$ and  $\delta$  be a
linear mapping from $\mathrm{Alg}\mathcal{L}$ into $B(X)$. Then
$\delta$ satisfies WJD property if and only if $\delta$ is a
generalized derivation and
$\delta(I)A\in(\mathrm{Alg}\mathcal{L})^\prime$ for every $A\in
\mathrm{Alg}\mathcal{L}$. In particular, if $\delta(I)=0$, then
$\delta$ satisfies WJD property if and only if $\delta$ is a
derivation.
\end{theorem}

\begin{corollary}\label{3003}
Let $\mathcal{L}$ be as in Example 1.1. Then
$\delta:\mathrm{Alg}\mathcal{L}\rightarrow B(H)$ satisfies WJD
property if and only if $\delta$ is a derivation.
\end{corollary}

\begin{proof}[\bf Proof]
By Theorem \ref{3002}, we only need to show that if $\delta$
satisfies WJD property, then $\delta(I)=0$. Let $n\geq2.$ By
\cite[Lemma 3.2]{C. Hou}, we have $(P_n)_-\ngeq P_n$. Hence there
exist
 $z_n\in P_n$ and $g_n\in (P_n)_-^\perp$ such that $g_n(z_n)=1$.
Also, there exists $y_n\in P_n$ such that $y_n$ and $z_n$ are linearly independent. Since
 $\delta$ satisfies WJD property,
we have $\delta(I)A\in\mathcal{A}^\prime$ for every $A\in \mathcal{A}$, which implies that
there exists some scalar $\lambda_n$
such that $\delta(I)x=\lambda_nx$ for every $x\in P_n$ and $\delta(I)(z_n\otimes g_n)
(y_n\otimes g_n)=\delta(I)(y_n\otimes g_n)(z_k\otimes g_n)$.
That is $\lambda_{n}g_n(y_n)z_n=\lambda_{n}y_n.$ The independence of $y_n$ and $z_n$
gives $\lambda_{n}=0$ and $\delta(I)x=0$ for every $x\in P_n$.
Since $\vee\{P_n\in\mathcal{L}:n=2,3,\cdots\}=H,$ it follows that $\delta(I)=0.$ The proof is complete.
\end{proof}

\begin{corollary}\label{3004}
Let $\mathcal{L}$ be a subspace lattice on  $H$ with $dimH\geq2$
such that $\vee\{L:L\in \mathcal{P}_\mathcal{L}\}=H$ or
$\wedge\{L_-:L\in \mathcal{P}_\mathcal{L}\}=(0)$. If $\mathcal{L}$
has a non-trivial comparable element, then
$\delta:\mathrm{Alg}\mathcal{L}\rightarrow B(H)$ satisfies WJD
property if and only if $\delta$ is a derivation.
\end{corollary}
\begin{proof}[\bf Proof]
According to Theorem \ref{3002}, we only need to show that if $\delta$
satisfies WJD property, then $\delta(I)=0$. By \cite[Proposition
2.9]{J. Li com}, we have
$(\mathrm{Alg}\mathcal{L})^\prime=\mathbb{C}I.$
 Hence by Theorem \ref{3002}, we have $\delta(I)=\lambda I$ and
$\delta(I)A=\mu_A I$ for every $A\in \mathrm{Alg}\mathcal{L}$ (where $\lambda,\mu_A\in \mathbb{C}$).
We claim that $\lambda=0.$ Suppose that $\lambda\neq 0,$
then every operator in $\mathrm{Alg}\mathcal{L}$ is a scalar multiple of the identity $I$. That is,
 for every $A\in \mathrm{Alg}\mathcal{L}$, the range of
$A$ is $H$ or $0$. However, Since $\mathrm{Alg}\mathcal{L}$ contains
a rank one operator, it is impossible.  Hence $\delta(I)=0.$
\end{proof}

By Corollary \ref{3004}, we can easily show the following result.

\begin{corollary}\label{3336}
Let
$\mathcal{L}$ be a subspace lattice on $H$
with $dimH\geq2$  such that
 $H_-\neq H$ or $(0)_+\neq(0)$.
Then $\delta:\mathrm{Alg}\mathcal{L}\rightarrow B(H)$ satisfies WJD
property if and only if $\delta$ is a derivation.
\end{corollary}

{\flushleft{\textbf{Remark.}} It follows from Theorems \ref{2004},
\ref{2005}, \ref{3002} and \ref{3333} that every  linear mapping
satisfying WJD property from a $\mathcal{P}$-subspace lattice
algebra into $B(X)$ is Jordan derivable  at zero point.} {But the
converse is not true. For example, let $\mathcal{T}_2(\mathbb{C})$
be the algebra of all $2\times2$ upper triangular matrices over the
complex field $\mathbb{C}$.  Define a linear mapping
$\delta:\mathcal{T}_2(\mathbb{C}) \rightarrow
\mathcal{T}_2(\mathbb{C})$ according to
$$\delta(\left(
          \begin{array}{cc}
            x_{11} & x_{12} \\
            0 & x_{22} \\
          \end{array}
        \right))=\left(
                  \begin{array}{cc}
                    x_{11} & x_{11}-x_{22}+x_{12} \\
                    0 & x_{22} \\
                  \end{array}
                \right),
$$
for all $x_{ij}\in \mathbb{C}$,$(1\leq i\leq j\leq 2)$. It is easy
to show that $\delta$ is a generalized derivation and
$\delta(I)=I\in(\mathcal{T}_2(\mathbb{C}))^\prime$, that is,
$\delta$ is Jordan derivable at zero point. However, it follows from
Corollary \ref{3336} that $\delta$ does not satisfy WJD property. }

\section{ Derivable mappings at zero point and local generalized derivations}

\

Let $\mathcal{A}$ be a unital algebra, $\mathcal{M}$ be an
$\mathcal{A}$-bimodule  and
 $\mathcal{T}$ be an ideal of $\mathcal{A}$.
We say that $\mathcal{T}$ is a \emph{left} (respectively \emph{right}) separating set of $\mathcal{M}$ if
 for every $m$ in $\mathcal{M}$, $m\mathcal{T}=\{0\}$ implies $m=0$
(respectively $\mathcal{T}m=\{0\}$ implies $m=0$). $\mathcal{T}$ is called
a \emph{separating set} of $\mathcal{M}$ if $\mathcal{T}$ is a left
separating set and a right separating set of $\mathcal{M}$. The
following result is obvious.

\begin{lemma}\label{4001}
Suppose that $\mathcal{L}$ is a subspace lattice on $X$ such that $\vee\{L:L\in \mathcal{P}_\mathcal{L}\}=X$
$($respectively $\wedge\{L_-:L\in \mathcal{P}_\mathcal{L}\}=(0)$$)$.
Then the ideal $\mathcal{T}=\mathrm{span}\{x\otimes f: x\in E, f\in E_-^\perp, E\in \mathcal{P}_\mathcal{L}\}$ of
$\mathrm{Alg}\mathcal{L}$ is a left $($respectively right$)$ separating set of $B(X)$.
\end{lemma}
%\begin{proof}[\bf Proof]
%Suppose that $\vee\{L:L\in \mathcal{P}_\mathcal{L}\}=X$. Let $A\in B(X)$ and $A\mathcal{T}=0.$
%For all $E\in\mathcal{P}_\mathcal{L}$, $x\in E$ and $f\in E_-^\perp$, we have $Ax\otimes f=0$.
%This means $Ax=0$ for all $x\in E$. Since $\vee\{L:L\in \mathcal{P}_\mathcal{L}\}=X$, we obtain $A=0$.
%Similarly, we can show that if $\wedge\{L_-:L\in \mathcal{P}_\mathcal{L}\}=(0)$, then $\mathcal{T}$ is
%a right separating set of $B(X)$.
%\end{proof}

By Lemmas \ref{2002} and \ref{4001}, we have the following result.

\begin{theorem}\label{4002}
Let $\mathcal{L}$ be a subspace lattice on $X$ such that $\vee\{L:L\in \mathcal{P}_\mathcal{L}\}=X$
or $\wedge\{L_-:L\in \mathcal{P}_\mathcal{L}\}=(0)$ and  $\delta$ be a linear mapping
from $\mathrm{Alg}\mathcal{L}$ into $B(X)$.
Then $\delta$ is derivable  at zero point if and only if $\delta$ is a  generalized
derivation and $\delta(I)\in (\mathrm{Alg}\mathcal{L})^\prime$.
In particular, if $\delta(I)=0,$ then $\delta$ is  derivable at zero point if and
only if $\delta$ is a derivation.
\end{theorem}

\begin{proof}[\bf Proof]
We will show that if $\mathcal{L}$ satisfies $\vee\{L:L\in \mathcal{P}_\mathcal{L}\}=X$
and $\delta$ is derivable at zero point,
then $\delta$ is a generalized derivation and $\delta(I)\in (\mathrm{Alg}\mathcal{L})^\prime$.
The proof for $\mathcal{L}$ with
$\wedge\{L_-:L\in \mathcal{P}_\mathcal{L}\}=(0)$ is similar.
By the proof of \cite[Lemma 3]{W. Jing S. Lu P. Li},
we may show that $\delta(AP)=\delta(A)P+A\delta(P)-A\delta(I)P$ and $\delta(I)P=P\delta(I),$
for every $A\in \mathrm{Alg}\mathcal{L}$ and every idempotent $P\in \mathrm{Alg}\mathcal{L}$.
With a proof similar to the proof of Claim 1 in Theorem \ref{2004},
we have $\delta(I)\in (\mathrm{Alg}\mathcal{L})^\prime$. Now for all
$A,B\in\mathrm{Alg}\mathcal{L}$ and $T\in\mathcal{T}$, we have
\begin{eqnarray*}
\delta(ABT)=\delta(AB)T+AB\delta(T)-AB\delta(I)T
\end{eqnarray*}
and
\begin{eqnarray*}
\delta(ABT)&=&\delta(A)BT+A\delta(BT)-A\delta(I)BT\\
&=&\delta(A)BT+A\delta(B)T+AB\delta(T)\\
&&-AB\delta(I)T-A\delta(I)BT.
\end{eqnarray*}
It follows that $\delta(AB)T=\delta(A)BT+A\delta(B)T-A\delta(I)BT.$ Since $\mathcal{T}$ is a left
 separating set of $B(X)$,
we obtain $\delta(AB)=\delta(A)B+A\delta(B)T-A\delta(I)B$ for all
$A,B\in\mathrm{Alg}\mathcal{L}$. That is, $\delta$ is a generalized
derivation. The proof is complete.
\end{proof}

Recall that a linear mapping $\delta$ from $\mathcal{A}$ into $\mathcal{M}$ is a
\emph{left} (respectively \emph{right}) \emph{multiplier}
if $\delta(AB)=\delta(A)B$ (respectively $\delta(AB)=A\delta(B)$) for all $A,B\in \mathcal{A}$; $\delta$ is
a \emph{local generalized derivation}
if for every $A\in \mathcal{A}$ there is a generalized derivation $\delta_A:\mathcal{A}\rightarrow \mathcal{M}$
(depending on $A$) such that $\delta(A)=\delta_A(A).$ In the following we give some applications
of Lemmas \ref{2002} and \ref{4001}.
The proofs of the results are similar to the proof of Theorem \ref{4002}, and we leave them to readers.

\begin{theorem}\label{4003}
Suppose that $\mathcal{L}$ is a subspace lattice on  $X$ such that $\vee\{L:L\in \mathcal{P}_\mathcal{L}\}=X$
$($respectively $\wedge\{L_-:L\in \mathcal{P}_\mathcal{L}\}=(0)$$)$ and $\delta$ is a linear
mapping from $\mathrm{Alg}\mathcal{L}$ into $B(X)$.
Then $\delta$ has the following properties:

$($a$)$ if $\delta(AB)=\delta(A)B$ for any $A,B\in \mathrm{Alg}\mathcal{L}$ with $AB=0$,
then $\delta$ is a left multiplier $($respectively
if $\delta(AB)=A\delta(B)$ for any $A,B\in \mathrm{Alg}\mathcal{L}$ with $AB=0$,
then $\delta$ is a right multiplier$)$;

$($b$)$ if $\delta(AB)=\delta(A)B+\delta(B)A$ for any $A,B\in \mathrm{Alg}\mathcal{L}$ with $AB=0$ and $\delta(I)=0$,
 then $\delta\equiv 0$ $($respectively if $\delta(AB)=A\delta(B)+B\delta(A)$ for any
 $A,B\in \mathrm{Alg}\mathcal{L}$ with $AB=0$ and $\delta(I)=0$,
 then $\delta\equiv 0$$)$;

$($c$)$ if $\delta(A^2)=2\delta(A)A$ for all $A\in \mathrm{Alg}\mathcal{L}$, then $\delta\equiv 0$
$($respectively if $\delta(A^2)=2A\delta(A)$ for all $A\in \mathrm{Alg}\mathcal{L}$, then $\delta\equiv 0$$)$.

\end{theorem}

Combining Theorem \ref{4003} (a) and \cite[Proposition 1.1]{J. Li Z. Pan}, we have

\begin{corollary}\label{4004}
Suppose that $\mathcal{L}$ is a subspace lattice on  $X$ such that $\vee\{L:L\in
 \mathcal{P}_\mathcal{L}\}=X$
and $\wedge\{L_-:L\in \mathcal{P}_\mathcal{L}\}=(0)$ and $\delta$ is a linear mapping
from $\mathrm{Alg}\mathcal{L}$ into $B(X)$.
Then the following are equivalent.

$($a$)$ $\delta$ is a generalized derivation.

$($b$)$ $\delta$ is a local generalized derivation.

$($c$)$ $A\delta(B)C=0$, whenever $A, B, C\in \mathrm{Alg}\mathcal{L}$ such that $AB=BC=0.$
\end{corollary}

Combining Lemmas \ref{2002}, \ref{4001} and \cite[Theorem 2.8]{J. Li Z. Pan J. Zhou}, we also have

\begin{theorem}\label{4005}
Let $\mathcal{L}$ be a subspace lattice on  $X$ such that $\vee\{L:L\in \mathcal{P}_\mathcal{L}\}=X$
and $\wedge\{L_-:L\in \mathcal{P}_\mathcal{L}\}=(0)$. If $h$ is a bijective linear mapping
from $\mathrm{Alg}\mathcal{L}$
onto a unital algebra satisfying $h(A)h(B)h(C)=0$ for all $A,B,C\in \mathrm{Alg}\mathcal{L}$
with $AB=BC=0$ and $\delta(I)=I,$ then $h$ is an isomorphism.
\end{theorem}

\

\end{document}